\documentclass[a4paper,12pt]{article}

\usepackage[latin1]{inputenc}
\usepackage[english]{babel}

\usepackage{amssymb}
\usepackage{amsmath}
\usepackage{amsthm}
\usepackage{amscd}
\usepackage{mathtools}
\usepackage{enumitem}

\usepackage{textcomp}
\usepackage{layout}
\usepackage{color}

\renewcommand{\Pr}{\mathbb P} 

\newcommand{\Poi}{\mathbf P} 

\newcommand{\Poy}{\mathsf S} 
\newcommand{\Rp}{\mathsf R} 

\newcommand{\Pp}{\mathsf P} 
\newcommand{\GP}{\mathsf D} 

\newcommand{\B}{\mathcal{B}} 
\newcommand{\Bbd}{\mathcal{B}_{0}} 

\newcommand{\Fsig}{\mathcal F}
\newcommand{\Gsig}{\mathcal G}

\newcommand{\MX}{\mathcal M(X)}

\newcommand{\MpmX}{\mathcal M^{\cdot\cdot}(X)}

\newcommand{\N}{\mathbb N}

\newcommand{\R}{\mathbb R}

\renewcommand{\d}{\mathrm{d}}

\renewcommand{\exp}{\operatorname{exp}}
\newcommand{\e}{\operatorname{e}}

\newcommand{\ccdot}{\,\cdot\,}

\theoremstyle{definition}
\newtheorem{defdefinition}{Definition}[section]

\theoremstyle{plain}
\newtheorem{defsatz}[defdefinition]{Theorem}
\newtheorem{defsatzdef}[defdefinition]{Theorem and Definition}
\newtheorem{defprop}[defdefinition]{Proposition}
\newtheorem{deflemma}[defdefinition]{Lemma}
\newtheorem{deffolgerung}[defdefinition]{Corollary}
\newtheorem{defbeispiel}[defdefinition]{Example}

\theoremstyle{remark}
\newtheorem{defbemerkung}[defdefinition]{Remark}

\newcommand{\satz}[1]{\begin{defsatz}#1\end{defsatz}}

\newcommand{\bem}[1]{\begin{defbemerkung}#1\end{defbemerkung}}
\newcommand{\lemma}[1]{\begin{deflemma}#1\end{deflemma}}
\newcommand{\lemman}[2]{\begin{deflemma}[#1]#2\end{deflemma}}

\newcommand{\korollarn}[2]{\begin{deffolgerung}[#1]#2\end{deffolgerung}}
\newcommand{\beispiel}[1]{\begin{defbeispiel}#1\end{defbeispiel}}

\newcommand{\F}{{\cal F}}
\newcommand{\A}{{\cal A}}

\title{Exit spaces for Cox processes and the P\'olya sum process}

\author{Mathias Rafler\footnote{email {\sf rafler@ma.tum.de}}}

\begin{document}

\maketitle

\begin{abstract}
For Cox processes we construct a Markov process with increasing paths to couple the condensations of the Cox process in a monotone way. A similar procedure procedure yields an analogue Markov process for the P\'olya sum process. Moreover, we identify the exit spaces of these Markov processes and identify them firstly as mixtures of certain extremal processes, i.e. as a process in a random environment, and secondly as Gibbs processes.
\end{abstract}

\section{Introduction}

A classical charcterization of Cox processes due to Mecke is one via condensations of point processes~\cite{jM11}. Such condensations are inverse operations of thinnings, and while a thinning always exists, condensation may not. Cox processes are according to Mecke's characterization exactly those point processes for which all condensations exist.

The first aim is to construct a Markov process on the state space of locally finite point measures with the following properties: Firstly, the sample paths are increasing; and secondly, the one-dimensional distributions of that Markov process are the condensations of the given Cox process. That way we obtain a coupling of Cox process and its condensations such that to realized point configurations points are only added.

The second basic aim is to construct the exit space of this Markov process yielding a representation of this `Cox Markov process' as a mixture of certain `Poisson Markov processes'. The directing measure turns out to be the one of the Cox process. These `Poisson Markov processes' are the Markov processes conditioned on the asymptotic $\sigma$-field. Along this way we get that the `Cox Markov processes' are Gibbs processes for some local specification in the sense of Preston~\cite{cP79}, and that the extremal states are exacly the `Poisson Markov processes'.

A particular Cox process is the P\'olya sum process, which is the Papangelou process for the P\'olya sum kernel. We present a similar construction of the Markov process corresponding to the P\'olya sum process with a different parametrization, since due to strong connections to the chinese restaurant process, this Markov process is of further interest. However, the presentation here restricts to the construction of the Martin-Dynkin boundary, the relations to the chinese restaurant process will be explored in~\cite{mR13scrp}.

Basically this presentation consists of two main parts, the discussion of Cox processes in section~\ref{sect:cox} and the discussion of the P\'olya sum process in section~\ref{sect:polya}. In each part we firstly discuss thinnings and splittings, then construct the corresponding Markov process and finally identify their exit spaces and identify these Markov processes as Gibbs processes. 

\section{Random measures and point processes}

Throughout the paper let $X$ be a Polish space with Borel $\sigma$-field $\B$ and a ring of bounded Borel sets $\Bbd$. By $F$ we denote the set of continuous functions from $X$ to $\R$, by $F_+\subset F$, $F_b\subset F$ the sets of non-negative, continuous functions and the set of continuous functions with bounded support, respectively. $F_{+,b}=F_+\cap F_b$.

When equipped with the vague topology, the space of locally finite measures $\MX$ is Polish as well as its closed subset of locally finite point measures $\MpmX$. Its Borel $\sigma$-field $\Gsig$ is generated by the evaluation mappings $\zeta_B:\MX\to\R_+$, $\zeta_B\mu=\mu(B)$. Let $\Gsig_B=\sigma(\zeta_{B'}:B\in\B, B\subset B)$. The elements $\MX$ are partially ordered via $\nu\leq\mu$ iff this inequality holds for all non-negative test functions on $X$. For point processes this means that the point masses of $\nu$ are dominated by the ones of $\mu$ or, equivalently, that $\mu-\nu$ is still a point measure.

A random measure is a probability measure on $\MX$, a point process is a random measure which is concentrated on $\MpmX$. For a random measure $\Rp$ the Campbell measure is defined as
\begin{equation*}
	C_\Rp(h)=\iint h(x,\nu)\nu(\d x)\Rp(\d\nu),
\end{equation*}
where $h:X\times\MX\to\R$ is measurable and non-negative or integrable. Point processes $\Pp$ satisfying some absolute continuity condition admit a particular disintegration of their Campbell measure,
\begin{equation*}
	C_\Pp(h)=\iint h(x,\mu+\delta_x)\pi(\mu,\d x)\Pp(\d\nu),
\end{equation*}
see e.g.~\cite{MWM79,NZ79,oK83,MKM78}. In this case, $\pi$ is called Papangelou kernel and $\Pp$ Papangelou process. Well known is the Poisson process $\Poi_\rho$ given by the kernel $\pi(\mu,\ccdot)=\rho$, $\rho\in\MX$. Of further special interest will be the P\'olya sum process $\Poy_{z,\rho}$ and the P\'olya difference process $\GP_{z,\rho}$ which are characterized by the sum kernel $\pi(\mu,\ccdot)=z(\rho+\mu)$, $z\in(0,1)$, $\rho\in\MX$, and the difference kernel $\pi(\mu,\ccdot)=z(\rho-\mu)$, $z>0$, $\rho\in\MpmX$, $\mu\leq\rho$. The latter processes were introduced and constructed in~\cite{hZ09,NZ11}.

Cox processes are Poisson processes with a random intensity measure. They obey a condenstion property due to Mecke~\cite{jM11}: A point process $\Pp$ is a Cox process if and only if for all $q\in(0,1]$ there exists a point process $P_q$ such that $\Pp$ is the $q$-thinning $\Gamma_q(P_q)$ of $P_q$, that is points of a realized configuration are deleted independently with probability $1-q$ and survive with probability $q$. Moreover, $P_q$ is unique and the directing measure of $\Pp$ is the weak limit of $P_q(q\ccdot)$ as $q\to 0$, i.e. it is obtained by increasing the number of points but decreasing their weights.

\section{Cox Markov processes and their exit spaces \label{sect:cox}}
\subsection{The Markov process}

We aim at constructing an increasing Markov process $Y=(Y_t)_{t\in[0,1)}$ with values in $\MpmX$ in the following way: Define $\Omega=\{\omega:[0,1)\to\MpmX: \omega\text{ cadlag, }\omega_s\leq\omega_t, s\leq t\}$ and equip $\Omega$ with the Skorohod topology. Define $Y_t(\omega)=\omega_t$. Moreover define two filtrations by
\begin{equation*}
	\Fsig_t=\sigma\bigl(Y_s:s\leq t\bigr),\qquad \Fsig^t=\sigma\bigl(Y_s:s\geq t\bigr).
\end{equation*}

In the following we fix a Cox process $\Pp$ and thus obtain according to Mecke's Theorem a family of point processes $(P_q)_{0<q\leq 1}$ such that $\Pp=\Gamma_q(P_q)$. Note that without the Cox property the following construction remains valid if there exists some $q_0>0$ such that $P_q$ exists for $q>q_0$ and everything is restricted to that smaller domain. Observe that in general a thinned process is a doubly stochastic P\'olya difference process, see also~\cite{NZ11}.

\lemman{Thinning}{ \label{thm:cox:thinning}
Let $0<q<q'<1$, then $P_{q'}=\Gamma_p(P_q)$ for $p=\tfrac{q}{q'}$.
}

\begin{proof}
Since
\begin{equation*}
	\Gamma_{q'}\bigl(\Gamma_p(P_q)\bigr) = \iint \GP_{\frac{q'}{1-q'},\mu} \GP_{\frac{p}{1-p},\nu}(\d\mu) P_q(\d\nu)
\end{equation*}
and 
\begin{equation*}
	\int \GP_{\frac{q'}{1-q'},\mu} \GP_{\frac{p}{1-p},\nu}(\d\mu)
		=\Gamma_{q'}(\GP_{\frac{p}{1-p},\nu})=\GP_{\frac{pq'}{1-p+p(1-q')},\nu}
\end{equation*}
by~\cite[Ex 6.3.1]{bN12}, the parameter of the thinning is obtained as $p=\tfrac{q}{q'}$.
\end{proof}

Following~\cite{bN12}, we define the splitting operation as
\begin{equation*}
	S_p\bigl(P\bigr)(h)=\iint h(\nu,\mu-\nu) \Gamma_p\bigl(\mu\bigr)(\d\nu)P(\d\mu).
\end{equation*}
By Lemma~\ref{thm:cox:thinning} we get the marginal distribution
\begin{equation*}
	S_p\bigl(P_q\bigr)(g\otimes 1)=\Gamma_p\bigl(P_q\bigr)(g),
\end{equation*}
and thus $S_p\bigl(P_q\bigr)(\ccdot\times A)\ll\Gamma_p\bigl(P_q\bigr)$, which allows the disintegration of $S_p(P_q)$,
\begin{equation*}
	S_p\bigl(P_q\bigr)(h)=\iint h(\nu,\eta) \Upsilon_p^\nu\bigl(P_q\bigr)(\d\eta)\Gamma_p\bigl(P_q\bigr)(\d\nu),
\end{equation*}
i.e. ensures the existence of the so-called splitting kernel $\Upsilon_p^{\ccdot}(P_q)$. In particular for $h(\nu,\eta)=g(\nu+\eta)$,
\begin{equation*}
	P_q(g)=S_p\bigl(P_q\bigr)(h)=\iint g(\nu+\eta) \Upsilon_p^\nu\bigl(P_q\bigr)(\d\eta) \Gamma_p\bigl(P_q\bigr)(\d\nu).
\end{equation*}
By using $\ast$ to denote the operation on the rhs, we obtain

\lemman{Splitting}{ \label{thm:cox:splitting}
Let  $0<q<q'<1$, then $P_q=\Upsilon_p^{\ccdot}(P_q)\ast P_{q'}$ for $p=\tfrac{q}{q'}$.
}

The splitting kernel should be understood as the distribution of the increment to get from $P_{q'}$ to $P_q$. Thinning and splitting are in some sense dual operations, which will be explored after the construction of the Markov process $Y=(Y_t)_{0\leq t<1}$.

Let
\begin{align}
	Y_0&\sim\Pp \label{eq:cox:mc:initial}\\
	\Pr(Y_t-Y_s\in\ccdot|Y_s) &=\Upsilon_{\frac{1-t}{1-s}}^{Y_s}(P_{1-t}),\qquad 0\leq s<t<1. \label{eq:cox:mc:transition}
\end{align}
Since the state space $\MpmX$ is complete and separable, the finite-dimensional distributions of $Y$ are uniquely determined, and by choosing cadlag paths, the process is uniquely determined. Because of Lemma~\ref{thm:cox:splitting}, $Y_t$ is $P_{1-t}$-distributed for each $t$. Furthermore, the transition probabilities of $Y$ are of the form
\begin{equation*}
	p_{s,t}(\nu,\ccdot)=\Upsilon_{\frac{1-t}{1-s}}^{\nu}(P_{1-t})\ast\Delta_\nu,
\end{equation*}
where $\Delta_\nu$ is the point process which realizes the point configuration $\nu$ a.s.
\bem{
The random initial condition might be unpleasent and one may prefer an initial condition like $Y_{t_0}=0$. This may be obtained immediatly by constructing $Y_t$ for $t\in[-1,0)$ as a thinning of $\Pp$ such that the law of $Y_t$ is $\Gamma_{t+1}(\Pp)$ and $Y_{-1}=0$, which can be done in an increasing way, too. However, this is rather irrelevant for the construction of its exit space and we omit that procedure. A parametrization of the P\'olya sum process in terms of one of its parameters yields the non-random initial condition in a natural way.
}

The backward dynamics follows from Lemma~\ref{thm:cox:thinning} and is given by
\begin{equation*}
	p^\ast_{s,t}(\nu,\phi)=\Gamma_{\frac{1-t}{1-s}}\bigl(\nu\bigr)(\phi)=\GP_{\frac{1-t}{t-s},\nu}(\phi).
\end{equation*}

\korollarn{Duality}{
Let $0\leq s<t<1$ and $\phi$, $\psi$ non-negative, measurable functions. Then
\begin{equation*}
	\int p_{s,t}^\ast(\ccdot,\psi)\phi \,\d P_{1-t}= \int \psi p_{s,t}(\ccdot,\phi) \,\d P_{1-s}.
\end{equation*}
}

\begin{proof}
\begin{align*}
	\int p_{s,t}^\ast(\mu,\psi)\phi(\mu) P_{1-t}(\d\mu)
		&= \iint \phi(\nu)\phi(\mu) \Gamma_{\frac{1-t}{1-s}}\bigl(\mu\bigr)(\d\nu) P_{1-t}(\d\mu)\\
		&= \iint \phi(\nu)\phi(\nu+\eta) \Upsilon_{\frac{1-t}{1-s}}^\nu\bigl(P_{1-t}\bigr)(\d\eta) \Gamma_{\frac{1-t}{1-s}}\bigl(P_{1-t}\bigr)(\d\nu)\\
		&= \int \psi(\nu) p_{s,t}(\nu,\phi) \Gamma_{\frac{1-t}{1-s}}\bigl(P_{1-t}\bigr)(\d\nu)
\end{align*}
due to the disintegration
\end{proof}

The family $(p_{s,t})_{0\leq s<t<1}$ defines a family of contraction operators $T=(T_{s,t})_{0\leq s\leq t<1}$ via
\begin{equation} \label{eq:operatorfamily}
	T_{s,t}\phi=p_{s,t}(\ccdot,\phi)
\end{equation}
for integrable functions $\phi:\MpmX\to\R$. When restricted to continuous functions which are measurable wrt. $\Gsig_B$ for some bounded $B$, then some properties hold.
\lemma{
Let $0\leq s\leq t<1$ and assume that $\nu\mapsto\Upsilon^\nu_{\frac{1-t}{1-s}}(P_{1-t})$ is vaguely contiunuous. Then if $\phi$ is continuous, also $t\mapsto T_{s,t}\phi$ is fore all $s$.
}
This follows from a standard argument.
\lemma{
If $\phi$ is continuous and measurble wrt. some $\Gsig_B$, then for all $s$,
\begin{equation*}
	\lim_{t\searrow s}T_{s,t}\phi=\phi.
\end{equation*}
}

\begin{proof}
Let $\phi$ be a bounded and $\Fsig_B$-measurable function. Then
\begin{align*}
	p_{s,t}(\nu,\phi)-\phi(\nu)&=\int \phi(\nu+\mu)\Upsilon_{\frac{1-t}{1-s}}^{\nu}\bigl(P_{1-t}\bigr)(\d\mu) -\phi(\nu)\\
		&=\begin{multlined}[t]
				\left[\Upsilon_{\frac{1-t}{1-s}}^{\nu}\bigl(P_{1-t}\bigr)(\zeta_B=0)-1\right] \phi(\nu) \\
				+\int_{\{\zeta_B>0\}} \phi(\nu+\mu) \Upsilon_{\frac{1-t}{1-s}}^{\nu}\bigl(P_{1-t}\bigr)(\d\mu).
		\end{multlined}
\end{align*}
Since $P_{1-t}$ converges weakly to $P_{1-s}$ as $t\searrow s$, certainly the splitting kernel converges weakly for a.e. $\nu$ to $\Delta_0$
\end{proof}

For finite Cox processes we may even compute the generator of the Markov process $Y$. Helpful is in this case the representation of the splitting kernel in terms of reduced Palm kernels~\cite[Thm. 6.3.5]{bN12}.
\lemma{
Let $\phi$ be a bounded and continuous function and assume that for each finite configuration $\nu$ the reduced Palm kernels $P_{1-t,\nu}^!\to P_{1-s,\nu}^!$ weakly as $t\searrow s$. Then the generator of $Y$ is given by
\begin{equation*}
	\A_s\phi(\nu)= \frac{1}{(1-s)P_{1-s,\nu}^!(\zeta_X=0)}\int_{\{\zeta_X=1\}} \phi(\nu+\mu)-\phi(\nu) P_{1-s,\nu}^!(\d\mu).
\end{equation*}
}
In particular this means that at each time only one point is allowed to arrive.
\begin{proof}
Let $\phi$ be a bounded and measurable function. Then
\begin{align*}
	p_{s,t}(\nu,\phi)-\phi(\nu)&=\int \phi(\nu+\mu)\Upsilon_{\frac{1-t}{1-s}}^{\nu}\bigl(P_{1-t}\bigr)(\d\mu) -\phi(\nu)\\
		&=\begin{multlined}[t]
				\left[\Upsilon_{\frac{1-t}{1-s}}^{\nu}\bigl(P_{1-t}\bigr)(\zeta_X=0)-1\right] \phi(\nu) \\
				+\int_{\{\zeta_X=1\}} \phi(\nu+\mu) \Upsilon_{\frac{1-t}{1-s}}^{\nu}\bigl(P_{1-t}\bigr)(\d\mu). \\
				+\int_{\{\zeta_X>1\}} \phi(\nu+\mu) \Upsilon_{\frac{1-t}{1-s}}^{\nu}\bigl(P_{1-t}\bigr)(\d\mu).
		\end{multlined}
\end{align*}
Divide by $t-s$ and denote these summands by $I_0\cdot\phi(\nu)$, $I_1$ and $I_2$. In the finite case
\begin{equation*}
	\Upsilon_{\frac{1-t}{1-s}}^{\nu}\bigl(P_{1-t}\bigr)
		=\frac{ \left(\frac{t-s}{1-s}\right)^{\zeta_X} } {P_{1-t,\nu}^!\Bigl(\bigl(\frac{t-s}{1-s}\bigr)^{\zeta_X}\Bigr)} P_{1-t,\nu}^!.
\end{equation*}
By $Z_{s,t}$ denote the normalization constant. Observe that as $t\searrow s$,
\begin{equation*}
	Z_{s,t}=\int \left(\frac{t-s}{1-t}\right)^{\mu(X)} P_{1-t,\nu}^!(\d\mu)
		\to P_{1-s,\nu}^!(\zeta_X=0)=:Z_s.
\end{equation*}
Thus for the second summand we get
\begin{align*}
	I_1 
		&= \frac{1}{1-t}\cdot\frac{1}{Z_{s,t}} \int_{\zeta_x= 1} \phi(\nu+\mu)  P_{1-t,\nu}^!(\d\mu)\\
		&\to  \frac{1}{1-s}\cdot\frac{1}{Z_{s}} \int_{\zeta_x= 1} \phi(\nu+\mu)  P_{1-s,\nu}^!(\d\mu).
\end{align*}
Analogously the third summand $I_2$ vanishes. Finally,
\begin{align*}
	I_0
		&=\frac{1}{Z_{s,t}} \int \frac{1}{t-s} \left(1_{\{\zeta_X=0\}}(\mu) -\left(\frac{t-s}{1-t}\right)^{\mu(X)}\right) P_{1-t,\nu}^!(\d\mu) \\
		&=\frac{1}{Z_{s,t}} \int_{\{\zeta_X>0\}} -\frac{(t-s)^{\mu(X)-1}}{(1-t)^{\mu(X)}} P_{1-t,\nu}^!(\d\mu)  \\
		&\to \frac{1}{1-s}\cdot\frac{1}{Z_s}\cdot P_{1-s,\nu}^!(\zeta_X=1).
\end{align*}
Final rearrangements yield the result.
\end{proof}

\subsection{The Exit space}

By $\F^1$ denote the asymptotic $\sigma$-field
\begin{equation*}
	\Fsig^1=\bigcap_{t<1} \Fsig^t.
\end{equation*}
Then we get for any $t<1$ and $\Fsig_t$-measurable function $\Phi$
\begin{equation*}
	\Pr(\Phi|\Fsig^1)=\lim_{t\to 1}\Pr(\Phi|\Fsig^t).
\end{equation*}
For $\sigma(Y_t)$-measurable random variables $\Phi$ the conditional expectation on the rhs. is for $1>T>t$,
\begin{equation} \label{eq:cox:cond_distr}
	\Pr\bigl(\phi(Y_t)|\Fsig^T\bigr)=\int \phi(\mu) \GP_{\frac{1-T}{T-t},Y_T}(\d\mu),
\end{equation}
which admits the calculation of the limit $T\to 1$.

\lemma{ \label{thm:cox:limit}
Let $\omega\in\Omega$ such that $(1-t)\omega_t$ converges vagely to some $\nu\in\MX$. Then for each continuous $f:X\to\R_+$ with bounded support
\begin{equation*}
	\lim_{T\to 1} \Pr\bigl(\e^{-Y_t(f)}|\Fsig^T\bigr)(\omega) = \exp\left(-\frac{1}{1-t}\int 1-\e^{-f}\d\nu\right)
\end{equation*}
}

\begin{proof}
By equation~\eqref{eq:cox:cond_distr} we get for $f\in F_+$
\begin{align*}
	\Pr\bigl(\e^{-Y_t(f)}|\Fsig^T\bigr)
		&=\exp\left(\int\log \frac{1+\frac{1-T}{T-t}\e^{-f}}{1+\frac{1-T}{T-t}} \d Y_T \right) \\
		&=\exp\left(\int\log \frac{T-t+(1-T)\e^{-f}}{1-t} \d Y_T \right) \\
		&=\exp\left(\int\log\left(1+(1-T) \frac{\e^{-f}-1}{1-t}\right)^{\frac{1}{1-T}} \d(1-T) Y_T \right).
\end{align*}
As $T\to 1$, the argument of the logarithm converges to $\exp\bigl[-\tfrac{1}{1-t}\bigl( 1-\e^{-f}\bigr)\bigr]$ and hence
\begin{equation*}
	\Pr\bigl(\e^{-Y_t(f)}|\Fsig^T\bigr)\to \exp\left( -\tfrac{1}{1-t} \int 1-\e^{-f}\d Q\right)
\end{equation*}
on the set $\{\omega:(1-t)\omega_t\text{ converges vagely}\}$ and $Q_\omega$ denotes this limit.
\end{proof}

\bem{
Note that if the one-dimensional distributions are replaced by $n$-dimensional ones at times $t_1<\ldots<t_n$, then since the thinning is Markovian, the previous lemma carries over to that case:
\begin{equation*}
	\begin{multlined}[t]
			\lim_{T\to 1} \Pr\bigl(\e^{-Y_{t_1}(f_1)}\cdots\e^{-Y_{t_n}(f_n)}|\Fsig^T\bigr) = \int \e^{-Y_{t_1}(f_1)}\cdots\e^{-Y_{t_n}(f_n)} \GP_{\tfrac{1-t_2}{t_2-t_1},\mu_2}(\d\mu_1) \times \\
			\cdots\times \GP_{\tfrac{1-t_n}{t_n-t_{n-1}},\mu_n}(\d\mu_{n-1}) \Poi_{\frac{1}{1-t_n}\nu}(\d\mu_n).
		\end{multlined}
\end{equation*}
We denote this probability by $\Pp_\nu$, by $Q_\nu=\lim (1-t)\omega_t$, and moreover by $\Rp$ the distribution of $Q$.
}

\satz{ \label{thm:cox:dblystochrep}
Let $\Phi$ be non-negative, $\Fsig_t$-measurable for some $t\in[0,1)$, then there exists a probability measure $\Rp$ on $\MX$ such that
\begin{equation*}
	\Pr(\Phi)=\int \Pp_\nu(\Phi) \Rp(\d\nu).
\end{equation*}
Moreover, $\Rp$ is the directing measure of the original Cox process $\Pp$.
}

\begin{proof}
Since the law of $Y_T$ is $P_{1-T}$, $(1-T)Y_T$ converges vaguely $\Pr$-a.s, hence the convergence in Lemma~\ref{thm:cox:limit} holds almost surely. Moreover, the law $\Rp$ of this limit is the directing measure of $\Pp$ as well as $\Pr$.
\end{proof}

\bem{
From a statistical mechanical point of view, the family of stochastic kernels
\begin{equation*}
	\pi_T(\omega,\ccdot)=\Pr\bigl(\ccdot|\Fsig^T\bigr)(\omega),\qquad T\in[0,1),
\end{equation*}
is a local specification in the sense of Preston~\cite{cP79}, see also~\cite{hF75}. By Theorem~\ref{thm:cox:dblystochrep} and Lemma~\ref{thm:cox:thinning}, each Markov process constructed from some Cox process via~\eqref{eq:cox:mc:initial}~and~\eqref{eq:cox:mc:transition} is consistent with $(\pi_T)_T$ in the sense that $\Pr\pi_T=\Pr$, i.e. they are Gibbsian. Moreover, thanks to Meckes characterization of Cox processes, there are no further compatible Markov processes.
}

\satz{
The family of stochastic kernels $(\pi_T)_{0\leq T<1}$ is a local specification and the set of extremal points consists of $\Pp_\nu$, $\nu\in\MX$.
}
We conclude this part with an example taken from~\cite{bN12}. Details on existence can be found therein. This example also includes the one to be discussed in the next section.
\beispiel{ \label{ex:cox:permanental}
For given $z\in(0,1)$ and a projective family kernels $(B^x_m)_{m\geq 1}$, where for each $x\in X$ and $m\in\N$, $B^x_m$ is a kernel from $X$ to $X^m$, such that for each $m$ the measure $B^{x_1}_{m-1}(\d x_2,\ldots,\d x_m)\lambda(\d x_1)$ is cyclic invariant, let ${\mathfrak I}_z$ be the infinitely divisible point process with Levy mesure 
\begin{equation*}
	L_z(\phi)=\sum_{m\geq 1} \frac{z^m}{m} \int_{X^m}\phi(\delta_{x_1}+\ldots+\delta_{x_m}) B^{x_1}_{m-1}(\d x_2,\ldots,\d x_m)\lambda(\d x_1).
\end{equation*}
Then ${\mathfrak I}_z$ is Cox and $Y_t\sim {\mathfrak I}_{\frac{z}{z+(1-z)(1-t)}}$.
}

\section{The P\'olya sum process \label{sect:polya}}

We next develop these ideas for a particular Cox process, the P\'olya sum process along the same lines. Note that instead of parametrizing wrt. the condensation parameter, we use its parameter $z$. 

The following two lemmas can be found in~\cite{bN12}, but we give proofs via partial integration at the end. 
\lemman{Sampling from the P\'olya sum processes}{ \label{thm:thinnings:polyathinning}
Let $\Pp=\Gamma_q(\Poy_{z,\rho})=\int\GP_{\frac{q}{1-q},\mu}\Poy_{z,\rho}(\d\mu)$. Then $\Pp=\Poy_{\gamma,\rho}$, where $\gamma=\gamma(z,q)=\tfrac{zq}{1-z(1-q)}$.
}%

\lemman{Condensation of the P\'olya sum processes}{ \label{thm:splittings:polyasplitting}
Let $0\leq\gamma\leq z<1$ and 
\begin{equation*}
  \Pp(\phi)=\iint \phi(\mu+\nu) \Poy_{\frac{z-\gamma}{1-\gamma},\rho+\nu}(\d\mu) \Poy_{\gamma,\rho}(\d\nu).
\end{equation*}
Then $\Pp=\Poy_{z,\rho}$.
}%
In this particular case the Markov process is
\begin{align*}
	Y_0&=0\\
	\Pr(Y_t-Y_s\in\ccdot|Y_s) &= \Poy_{\frac{t-s}{1-s},\rho+Y_s},\qquad 0\leq s<t<1.
\end{align*}
Because of Lemma~\ref{thm:splittings:polyasplitting}, $Y_t$ is $\Poy_{t,\rho}$-distributed for each $t$. Furthermore, the transition probabilities of $Y$ are of the form
\begin{equation*}
	p_{s,t}(\nu,\ccdot)=\Poy_{\frac{t-s}{1-s},\rho+\nu}\ast\Delta_\nu.
\end{equation*}

To construct the exit space we need the backward dynamics $p^\ast$, but we already know that this is just done by independent thinning,
\begin{equation*}
	p^\ast_{s,t}(\nu,\ccdot)=\Gamma_{\frac{s(1-t)}{t(1-s)}}(\nu)
\end{equation*}
for $s<t$ (one starts at time $t$ with configuration $\nu$ and goes back to times $s$). By the sampling lemma this means that 
\begin{equation*}
	p^\ast_{s,t}(\nu,\phi)=\GP_{\frac{s(1-t)}{t-s},\nu}(\phi).
\end{equation*}

\lemma{
Let $\omega\in\Omega$ such that $(1-t)\omega_t$ converges vagely to some $\nu\in\MX$. Then for each continuous $f:X\to\R_+$ with bounded support
\begin{equation*}
	\lim_{T\to 1}\Pr\bigl(\e^{-Y_t(f)}|\Fsig^T\bigr)(\omega) 
		= \exp\left(-\frac{t}{1-t}\int 1-\e^{-f}\d\nu\right)
\end{equation*}
}
\begin{proof}
By the sampling lemma, $\Pr(\phi(Y_t)|\Fsig^T)=\int\phi(\mu) \GP_{\frac{t(1-T)}{T-t},Y_T}(\d\mu)$ for measurable $\phi$ and we get
\begin{align*}
	\Pr\bigl(\e^{-Y_t(f)}|\Fsig^T\bigr)
		&=\exp\left(\int\log \frac{1+\frac{t(1-T)}{T-t}\e^{-f}}{1+\frac{t(1-T)}{T-t}} \d Y_T \right) \\
		&=\exp\left(\int\log \frac{T-t+t(1-T)\e^{-f}}{T(1-t)} \d Y_T \right) \\
		&=\exp\left(\int\log\left(1+(1-T) \frac{t(\e^{-f}-1)}{T(1-t)}\right)^{\frac{1}{1-T}} \d(1-T) Y_T\right)\\
\end{align*}
As $T\to 1$, the argument of the logarithm converges to $\exp\bigl[-\tfrac{t}{1-t}\bigl( 1-\e^{-f}\bigr)\bigr]$ and hence
\begin{equation*}
	\Pr\bigl(\e^{-Y_t(f)}|\Fsig^T\bigr)\to \exp\left( -\tfrac{t}{1-t} \int 1-\e^{-f}\d Q\right)
\end{equation*}
on the set $\{\omega:(1-t)\omega_t\text{ converges vagely}\}$ and $Q$ denotes this limit.
\end{proof}


\bem{
The results for $\Phi(\omega)=\prod_{j=1}^{n}\e^{-\omega_{t_j}(f_j)}$ for $t_1<t_2<\cdots<t_n$, $n\in\N$ followdirectly from that lemma since for the limit only $t_n$ matters.
}

Thus define $Q_\omega=\lim_{t\to 1}(1-t)\omega_t$ as the vague limit in case of existence and 0 otherwise. Then $Q$ is $\Fsig^1$-measurable and
\begin{equation*}
	\Pr\bigl(\e^{-Y_t(f)}|\Fsig^1\bigr)(\omega)=\Poi_{Q_\omega}(\e^{-\zeta_f})
\end{equation*}
if the mentioned limit $Q_\omega$ exists, where $\Poi_{\nu}$ denotes the Poisson process with intensity measure $\nu$. But this limit exists a.s. and its distribution is the Poisson-Gamma process, hence

\satz{
For non-negative, $\Fsig_t$-measurable $\Phi$ for some $t\in[0,1)$, 
\begin{equation*}
	\Pr(\Phi)=\int \Pp_\nu(\Phi) \Rp(\d\nu)
\end{equation*}
where $\Rp$ is the Poisson-Gamma process the random measure $\Rp$ being the solution of the integration parts formula
\begin{equation*}
	C_\Rp(h)=\iint h(x,\mu+r\delta_x)\e^{-r}\d r\rho(\d x)\Rp(\d\mu).
\end{equation*}
}

\bem{
For fixed $B\in\Bbd$, the process of the marginals $\bigl(Y_t(B)\bigr)_t$ is a Markov process with values in $\N_0$ such that for each $t$, $Y_t(B)\sim NB\bigl(\rho(B),t\bigr)$. The transition from $Y_s(B)$ to $Y_t(B)$ means to add a $NB\bigl(\rho(B)+Y_s(B),\tfrac{t-s}{1-s}\bigr)$ distributed random variable, and backwards to take $Y_s(B)$ conditioned on $Y_t(B)$to be $B\bigl(Y_t(B),\frac{s(1-t)}{t(1-s)}\bigr)$ distributed.
}

\bem{
In a similar way further point processes may be analyzed
\begin{enumerate}
	\item The example~\ref{ex:cox:permanental} can be treated in the same way as the P\'olya sum process with the initial configuration being the empty configuration. However, an explicit representation of the splitting kernel is currently not available.
	\item For a similar Poisson Markov process, we have $p_{s,t}(\nu,\phi)=\int \phi(\nu+\mu)\Poi_{(t-s)\rho}(\d\mu)$, whereas the backwards dynamics is again given by the independent thinning $p^\ast_{t,\nu}(s,\phi)=\GP_{\frac{s}{t-s},\nu}(\phi)$, in both cases $0\leq s\leq t<\infty$. Then applying the limit procedure to 
		\begin{align*}
			\Pr\bigl(\e^{-Y_t(f)}|\Fsig^T\bigr) &=\exp\left(\int\log \frac{1+\frac{t}{T-t}\e^{-f}}{1+\frac{t}{T-t}} \d Y_T \right)\\
				& =\exp\left(\int\log\left(1+ \frac{t}{T}(\e^{-f}-1)\right)^T \d\frac{Y_T}{T} \right)
		\end{align*}
		yields that the law of $Y_t$ is a Poisson process with intensity measure $tQ$, where $Q=\lim_{T\to\infty}\tfrac{Y_T}{T}$. Since we started with a Poisson process, $Q=\rho$ a.s.

		\item The P\'olya sum process may be replaced by the P\'olya difference process: Then firstly the $q$-thinning is $\Gamma_q(\GP_{z,\rho})=\GP_{\frac{zq}{1+z(1-q)},\rho}$ and therefore $p^\ast_{s,t}(\nu,\ccdot)=\Gamma_{\frac{s(1+t)}{t(1+s)}}(\nu)=\GP_{\frac{s(1+t)}{t-s},\nu}$, secondly the splitting property turns into $\GP_{\frac{z-\gamma}{1+\gamma},\rho-\ccdot}\ast\GP_{\gamma,\rho}=\GP_{z,\rho}$ and $p_{s,t}(\nu,\ccdot)=\GP_{\frac{t-s}{1+s},\rho-\ccdot}\ast\Delta_\nu$, where now $0\leq s<t<\infty$. Carrying out the limit procedure, one obtains that $\Pr(\phi(Y_t)|\Fsig^\infty)=\lim_{T\to\infty}\GP_{\frac{t(1+T)}{T-t},Y_T}(\phi)=\lim_{T\to\infty}\Gamma_{\frac{t(1+T)}{T(1+t)}}\bigl(Y_T\bigr)(\phi)=\Gamma_{\frac{t}{1+t}}\bigl(Q\bigr)(\phi)$, where $Q=\lim_{T\to\infty}Y_T$ in case of existence. But with underlying $\GP_{z,\rho}$ this limit is $\rho$ a.s.
\end{enumerate}
}

\section{Scholion: Proofs of Sampling and condensation lemmas}


Proofs of these two lemmas can be found in Nehring's thesis~\cite{bN12}. Here we give proofs 
in terms of partial integration.
\begin{proof}[Proof of Lemma~\ref{thm:cox:thinning}]
Denote $\GP_\mu=\GP_{\frac{q}{1-q},\mu}$ and $u=\tfrac{q}{1-q}$ for short
\begin{align*}
	C^!_\Pp(h) &= \int C^!_{\GP_\mu}(h) \Poy_{z,\rho}(\d\mu)
		= u\iiint h(x,\nu)\bigl(\mu-\nu\bigr)(\d x) \GP_\mu(\d\nu) \Poy_{z,\rho}(\d\mu)\\
		&=\begin{multlined}[t]
				uz\iiint h(x,\nu)\GP_{\mu+\delta_x}(\d\nu) \bigl(\rho+\mu\bigr)(\d x) \Poy_{z,\rho}(\d\mu)\\
				-u\iiint h(x,\nu)\nu(\d x)\GP_{\mu}(\d\nu) \Poy_{z,\rho}(\d\mu).
			\end{multlined}
\end{align*}
Note that $\GP_{\mu+\delta_x}=\GP_\mu\ast\GP_{\delta_x}$ and $\GP_{\delta_x}(\phi)=(1-q)\phi(0)+q(\phi(\delta_x)$. Thus
\begin{align*}
	C^!_\Pp(h) &=\begin{multlined}[t]
				zq\iiint h(x,\nu)\GP_{\mu}(\d\nu) \bigl(\rho+\mu\bigr)(\d x) \Poy_{z,\rho}(\d\mu)\\
				+uzq\iiint h(x,\nu+\delta_x)\GP_{\mu}(\d\nu) \bigl(\rho+\mu\bigr)(\d x) \Poy_{z,\rho}(\d\mu)\\
				-u\iiint h(x,\nu)\nu(\d x)\GP_{\mu}(\d\nu) \Poy_{z,\rho}(\d\mu).
			\end{multlined}
\end{align*}
By introducing an alimentary null and applying the partial integration the first summand turns into
\begin{equation*}
	zq\iiint h(x,\nu) \bigl(\rho+\nu\bigr)(\d x) \GP_{\mu}(\d\nu) \Poy_{z,\rho}(\d\mu)
		+z(1-q) C^!_\Pp(h).
\end{equation*}
By the same procedure, the second summand turns into
\begin{equation*}
	uzq\iiint h(x,\nu+\delta_x) \bigl(\rho+\nu\bigr)(\d x) \GP_{\mu}(\d\nu) \Poy_{z,\rho}(\d\mu)
		+zq C_\Pp(h).
\end{equation*}
Therefore,
\begin{align*}
	\bigl(1-z(1-q)\bigr)C^!_\Pp(h) &=\begin{multlined}[t]
				zq\iint h(x,\nu) \bigl(\rho+\nu\bigr)(\d x) \Pp(\d\nu)\\
				+uzq\iint h(x,\nu+\delta_x) \bigl(\rho+\nu\bigr)(\d x) \Pp(\d\nu)\\
				-u\bigl(1-z(1-q)\bigr)C_\Pp(h).
			\end{multlined}
\end{align*}
Thus
\begin{equation*}
	C^!_\Pp(g) = \gamma \iint g(x,\nu) \bigl(\rho+\nu\bigr)(\d x) \Pp(\d\nu)
\end{equation*}
for $g(x,\mu)=h(x,\mu)+uh(x,\mu+\delta_x)$.
\end{proof}


\begin{proof}[Proof of Lemma~\ref{thm:cox:splitting}]
Write $u=\tfrac{z-\gamma}{1-\gamma}$ and $v=\tfrac{\gamma(1-z)}{1-\gamma}$, then $u=z-\gamma v$ and $1-u=v$. Let $h$ be integrable, then by partial integration
\begin{align*}
  C_\Pp(h) &= \iiint h(x,\mu+\nu) \bigl(\mu+\nu\bigr)(\d x) \Poy_{u,\rho+\nu}(\d\mu) \Poy_{\gamma,\rho}(\d\nu)\\
  	&=\begin{multlined}[t]
  			u\iiint h(x,\mu+\nu+\delta_x) \bigl(\rho+\nu+\mu\bigr)(\d x) \Poy_{u,\rho+\nu}(\d\mu) \Poy_{\gamma,\rho}(\d\nu)\\
  			+ \gamma \iiint h(x,\mu+\nu+\delta_x) \Poy_{u,\rho+\nu+\delta_x}(\d\mu) \bigl(\rho+\nu\bigr)(\d x) \Poy_{\gamma,\rho}(\d\nu).
  		\end{multlined}
\end{align*}
To transform the second summand, note that $\Poy_{u,\rho+\nu+\delta_x} =\Poy_{u,\rho+\nu} \ast \Poy_{u,\delta_x}$, and moreover
\begin{equation*}
	\Poy_{u,\delta_x}(\phi)= (1-u)\sum_{n\geq 0}u^n \phi(n\delta_x).
\end{equation*}
Thus
\begin{align*}
  C_\Pp(h) &=\begin{multlined}[t]
  			z\iiint h(x,\mu+\nu+\delta_x) \bigl(\rho+\nu+\mu\bigr)(\d x) \Poy_{u,\rho+\nu}(\d\mu) \Poy_{\gamma,\rho}(\d\nu)\\
  			-\gamma v\iiint h(x,\mu+\nu+\delta_x) \bigl(\rho+\nu+\mu\bigr)(\d x) \Poy_{u,\rho+\nu}(\d\mu) \Poy_{\gamma,\rho}(\d\nu)\\
  			+ \gamma v \sum_{n\geq 0} u^n\iiint h(x,\mu+\nu+\delta_x) \bigl(\rho+\nu\bigr)(\d x) \Poy_{u,\rho+\nu}(\d\mu) \Poy_{\gamma,\rho}(\d\nu).
  		\end{multlined}
\end{align*}
Finally observe that the last line kills the second.
\end{proof}




\begin{thebibliography}{99}

\bibitem{eD65} Dynkin, E.~B. (1965).\newblock Markov processes, vol.1. \newblock Springer.

\bibitem{EK05} Ethier S.~N. and Kurtz, T.~G. (2005). \newblock {Markov Processes: Characterization and Convergence}. \newblock John Wiley \& Sons.


\bibitem{EGW12} Evans, S.N., Gr\"ubel, R. and Wakolbinger, A. (2012). \newblock Trickle-down processes and their boundaries. \newblock {\em Electron. J. Probab.}, 17(1):1--58.

\bibitem{hF75} F\"ollmer, H. (1975). \newblock Phase transition and Martin Boundary. \newblock {\em Seminaire de probabilites (Strasbourg)}, 9:305--17.
  


\bibitem{oK83} Kallenberg, O. (1983). \newblock {Random measures}. \newblock Akademie-Verlag and Academic Press.

\bibitem{MKM78} Kerstan, J., Matthes, K. and Mecke, J. (1978). \newblock {Infinitely Divisible Point Processes}. \newblock John Wiley \& Sons.





\bibitem{MWM79} Matthes, K., Warmuth, W. and Mecke, J. (1979). \newblock Bemerkungen zu einer Arbeit von Nguyen Xuan Xanh and Hans Zessin. \newblock {\em Math. Machr.}, 88:117--27.


\bibitem{jM11}
Mecke, J. (2011).
\newblock Random measures.
\newblock Walter Warmuth Verlag.

\bibitem{bN12} Nehring, B. (2012). \newblock Point processes in statistical mechanics: A cluster expansion approach. \newblock {\em Thesis}.

\bibitem{NZ11} Nehring, B. and Zessin, H. (2011). \newblock {The Papangelou process. A concept for Gibbs, Fermi and Bose processes}. \newblock {\em J. Contemp. Math. Anal.} 46:326--37.

\bibitem{NZ79} Nguyen, X.~X. and Zessin, H. (1979). \newblock Integral and differential characterisation of the Gibbs process. \newblock {\em Math. Nachr.}, 88:105--15. 

\bibitem{cP79} Preston, C. (1979). \newblock Canonical and Microcanonical Gibbs States. \newblock {\em Z. Wahrscheinlichkeitstheorie verw. Gebiete}, 46:125--58.
  
\bibitem{mR13scrp} Rafler, M. (2013). \newblock A hydrodynamic limit and fluctuations for a Chinese restaurant-like process. \newblock {\em preprint}.

\bibitem{hZ09} Zessin, H. (2009). \newblock Der Papangelou Proze{\ss}. \newblock {\em J. Contemp. Math. Anal.} 44:61--72.

\end{thebibliography}
\end{document}